\documentclass[english]{smfart}

\renewcommand{\guillemotright}%
  {\nobreak\leavevmode
   \hbox{\fontencoding{U}\fontfamily{lasy}%
         \fontseries{m}\fontshape{n}\selectfont
   \kern+0.20em)\kern-0.20em)}}

\usepackage{amssymb}
\usepackage[mathscr]{eucal}	

\newtheorem{thm}{Theorem}
\newtheorem{lem}[thm]{Lemma}
\newtheorem*{maintheorem}{Main Theorem \textup(Rapoport's Conjecture\textup)}

\renewcommand{\theenumi}{\roman{enumi}}
\renewcommand{\labelenumi}{(\theenumi)}
\numberwithin{equation}{section}

\usepackage[arrow,matrix,tips,frame,curve,ps,dvips]{xy}
\SelectTips{cm}{}
\UseTips
\CompileMatrices

\usepackage{wrapfig}
\newcommand{\CC}{{\mathbb C}}
\newcommand{\RR}{{\mathbb R}}
\newcommand{\QQ}{{\mathbb Q}}
\newcommand{\EE}{{\mathbb E}}
\newcommand{\HH}{{\mathbb H}}
\newcommand{\VV}{{\mathbb V}}

\DeclareMathOperator{\codim}{codim}
\DeclareMathOperator{\Cone}{Cone}
\DeclareMathOperator{\rank}{rank}
\DeclareMathOperator{\Ker}{Ker}
\renewcommand{\Im}{\operatorname{Image}}
\DeclareMathOperator{\QQrank}{\QQ-rank}
\DeclareMathOperator{\RRrank}{\RR-rank}
\DeclareMathOperator{\mS}{SS}
\DeclareMathOperator{\emS}{\mS_{ess}}
\DeclareMathOperator{\Res}{Res}
\DeclareMathOperator{\Type}{Type}
\renewcommand{\l}{\ell}

\newcommand{\lsb}[1]{{}_{\scriptscriptstyle #1}}
\newcommand{\lsp}[1]{{}^{#1}}
\newcommand{\tildearrow}{\xrightarrow{\sim}}

\newcommand{\G}{\Gamma}
\DeclareMathOperator{\SL}{SL}
\newcommand{\Sympl}{\operatorname{Sp}}
\DeclareMathOperator{\GSp}{GSp}
\newcommand{\Pl}{\mathscr P}

\newcommand{\Xbar}{\overline{X}}
\newcommand{\Xhat}{\widehat{X}}
\newcommand{\back}{\backslash}
\newcommand{\Dstar}{D^*}
\newcommand{\Xstar}{X^*}

\newcommand{\ihat}{{\hat{\imath}}}
\newcommand{\jhat}{{\hat{\jmath}}}

\newcommand{\IC}{{\mathcal I\mathcal C}}
\newcommand{\Sheaf}{\mathcal S}
\newcommand{\WC}{\mathcal W \mathcal C}
\renewcommand{\L}{\mathscr L}
\newcommand{\M}{\mathcal M}

\newcommand{\sa}{{\mathfrak a}}
\newcommand{\h}{{\mathfrak h}}
\newcommand{\hb}{{\mathfrak b}}
\renewcommand{\k}{{\mathfrak k}}
\newcommand{\n}{{\mathfrak n}}
\newcommand{\levi}{{\mathfrak l}}
\newcommand{\p}{{\mathfrak p}}

\newcommand{\al}{\alpha}
\newcommand{\D}{\Delta}
\newcommand{\g}{\gamma}
\renewcommand{\u}{\mu}
\renewcommand{\r}{\rho}
\renewcommand{\t}{\tau}
\renewcommand{\o}{\omega}

\DeclareMathOperator{\Rep}{\mathfrak M\mathfrak o\mathfrak d}
\DeclareMathOperator{\IrrRep}{\mathfrak I\mathfrak r\mathfrak r}

\let\Der\Derived
\newcommand{\X}{\mathcal X}
\DeclareMathOperator{\Complex}{\mathbf C}
\DeclareMathOperator{\Graded}{\mathbf G\mathbf r}

\author{Leslie Saper} \address{Department of Mathematics\\ Duke
University\\ Box 90320\\ Durham, NC 27708\\U.S.A.}
\email{saper@math.duke.edu}
\urladdr{http://www.math.duke.edu/faculty/saper}
\title[$\mathscr L$-modules and Rapoport's Conjecture]{$\mathscr
L$-modules and the Conjecture of Rapoport and Goresky-MacPherson}
\alttitle{$\mathscr
L$-modules et la Conjecture de Rapoport et Goresky-MacPherson}

\begin{document}

\frontmatter

\begin{abstract}
Consider the middle perversity intersection cohomology groups of various
compactifications of a Hermitian locally symmetric space.  Rapoport and
independently Goresky and MacPherson have conjectured that these groups
coincide for the reductive Borel-Serre compactification and the
Baily-Borel-Satake compactification.  This paper describes the theory of
$\mathscr L$-modules and how it is used to solve the conjecture.  More
generally we consider a Satake compactification for which all real boundary
components are equal-rank.  Details will be given elsewhere
\cite{refnSaperLModules}.  As another application of $\mathscr L$-modules,
we prove a vanishing theorem for the ordinary cohomology of a locally
symmetric space.  This answers a question raised by Tilouine.
\end{abstract}

\begin{altabstract}
Consid\'erons les groupes de cohomologie d'intersection
(de perversit\'e interm\'ediare) de diverses compactifications d'un espace
localement hermitien sym\'etrique.  Rapoport et, ind\'ependamment, Goresky et
MacPherson ont conjectur\'e que ces groupes co\"\i ncident pour la
compactification de Borel-Serre r\'eductive et la compactification de
Baily-Borel-Satake.  Cet article d\'ecrit la th\'eorie des $\mathscr
L$-modules et la fa\c con dont elle peut s'employer pour r\'esoudre la
conjecture.  Plus g\'en\'eralement, nous tra\^itons une compactification de
Satake pour laquelle toutes les composantes r\'eelles \`a la fronti\`ere
sont de \guillemotleft rang \'egal\guillemotright.  Les d\'etails en seront
disponibles ailleurs \cite{refnSaperLModules}. Comme application
suppl\'ementaire de la th\'eorie des $\mathscr L$-modules, nous prouvons un
th\'eor\`eme d'annulation sur le groupe de cohomologie ordinaire d'un
espace localement sym\'etrique.  Ceci r\'epond \`a une question soulev\'ee
par Tilouine.
\end{altabstract}

\thanks{This research was supported in part by the National Science
Foundation under grants DMS-8957216, DMS-9100383, and DMS-9870162, and by a
grant from the Institut des Hautes \'Etudes Scientifiques.  Early stages of
this research were supported by The Duke Endowment, a Alfred P. Sloan
Research Fellowship, and the Maximilian-Bickhoff-Stiftung during a visit to
the Katholischen Universit\"at Eichst\"att.  The author wishes to thank all
these organizations for their hospitality and support.}

\thanks{The manuscript was prepared with the \AmS-\LaTeX\ macro system;
diagrams were generated with the \Xy-pic\ package.}

\subjclass{11F75, 14G35, 22E40, 22E45, 32S60, 55N33}

\keywords{Intersection cohomology, Shimura varieties, locally symmetric
spaces, compactifications}

\altkeywords{cohomologie d'intersection, vari\'et\'es de Shimura, espace
localement sym\'etrique, compactifications}

\maketitle


\mainmatter
\section{Introduction}
In a letter to Borel in 1986 Rapoport made a conjecture (independently
rediscovered by Goresky and MacPherson in 1988) regarding the equality
of the intersection cohomology of two compactifications of a locally
symmetric variety, the reductive Borel-Serre compactification and the
Baily-Borel compactification.  In this paper I describe the
conjecture, introduce the theory of $\L$-modules which was developed
to attack the conjecture, and explain the solution of the conjecture.
The theory of $\L$-modules actually applies to the study of many other
types of cohomology.  As a simple illustration, I will answer at the
end of this paper a question raised during the semester by Tilouine
regarding the vanishing of the ordinary cohomology of a locally
symmetric variety below the middle degree.  Except in this final
section, proofs are omitted; the details will appear in
\cite{refnSaperLModules}.

This paper is an expanded version of lectures I gave during the
Automorphic Forms Semester (Spring 2000) at the Centre \'Emile Borel
in Paris; I would like to thank the organizers for inviting me and
providing a stimulating environment.  During this research I benefited
from discussions with numerous people whom I would like thank, in
particular A. Borel, R. Bryant, M. Goresky, R. Hain, G. Harder,
J.-P. Labesse, J. Tilouine, M. Rapoport, J. Rohlfs, J. Schwermer, and
N. Wallach.

\section{Compactifications}
We consider a connected reductive algebraic group $G$ defined over
$\QQ$ and its associated symmetric space $D = G(\RR)/K A_G$, where $K$
is a maximal compact subgroup of $G(\RR)$ and $A_G$ is the identity
component of the $\RR$-points of a maximal $\QQ$-split torus in the
center of $G$.  Let $\G\subset G(\QQ)$ be an arithmetic subgroup which
for simplicity here we assume to be neat.  (Any arithmetic subgroup
has a neat subgroup of finite index; the neatness hypothesis ensures
that all arithmetic quotients in what follows will be smooth as
opposed to $V$-manifolds or orbifolds.) The locally symmetric space
$X= \G\back D$ is in general not compact and we are interested in
three compactifications (see Figure~\ref{figCompactifications}),
belonging respectively to the topological, differential geometric, and
(if $D$ is Hermitian symmetric) complex analytic categories.

\begin{figure}
\entrymodifiers={++!}
\begin{equation*}
\xymatrix{
*+!!/d.8ex/{\txt{Borel-Serre}} & {\Xbar} \ar[d] \ar@{}[r]^-*!/d.5ex/{=} &
{\smash[b]{\coprod_{P\in\Pl} Y_P,}}
\ar@<3ex>[d]^{\text{\footnotesize Collapse $\G_{N_P}\back N_P(\RR)$ fibers}}
& {Y_P = \G_P\back P(\RR)/K_PA_P} \\ 
*+!!/d.8ex/{\txt{Reductive\\Borel-Serre}} & {\Xhat} \ar[d]_\pi \ar@{}[r]^-*!/d.5ex/{=} &
{\smash[b]{\coprod_{P\in\Pl} X_P,\vphantom{\Xhat}}}
\ar@<3ex>[d]^{\text{\footnotesize Project $X_P=
X_{P,\l}\times X_{P,h} \to X_{P,h} = X_{R,h}$}} & {X_P = \G_{L_P}\back L_P(\RR)/K_PA_P}
\\
*+!!/d.8ex/{\txt{Baily-Borel\\Satake}} & {\Xstar} \ar@{}[r]^-*!/d.5ex/{=} &
{\smash[t]{\coprod_{R\in\Pl_1} F_R,}\vphantom{\Xstar}} & {F_R = X_{R,h}}
}
\end{equation*}
\caption{}
\label{figCompactifications}
\end{figure}

Let $\Pl$ (resp. $\Pl_1$) denote the partially ordered set of
$\G$-conjugacy classes of parabolic (resp. maximal parabolic)
$\QQ$-subgroups of $G$.  For $P\in \Pl$, let $L_P$ denote the Levi quotient
$P/N_P$, where $N_P$ is the unipotent radical of $P$.  (When it is
convenient we will identify $L_P$ with a subgroup of $P$ via an appropriate
lift.)  The {\itshape Borel-Serre compactification\/} \cite{refnBorelSerre} has
strata $Y_P = \Gamma_P\back P(\RR) /K_P A_P$ indexed by $P\in\Pl$ (for
$P=G$ we simply have $Y_G=X$).  Here $\Gamma_P = \Gamma\cap P$, $K_P=K\cap
P$, and $A_P$ is the identity component of the $\RR$-points of a maximal
$\QQ$-split torus in the center of $L_P$.  The Borel-Serre compactification
$\Xbar$ is a manifold with corners, homotopically equivalent with $X$
itself.

The arithmetic subgroup $\G$ induces arithmetic subgroups $\Gamma_{N_P} =
\Gamma \cap N_P$ in $N_P$ and $\Gamma_{L_P} = \G_P/\Gamma_{N_P}$ in $L_P$.
Let $D_P = L_P(\RR)/K_PA_P$ be the symmetric space associated to $L_P$ and
let $X_P= \G_{L_P}\back D_P$ be its arithmetic quotient.  Each stratum of
$\Xbar$ admits a fibration $Y_P \to X_P $ with fibers being compact
nilmanifolds $\Gamma_{N_P}\back N_P(\RR)$.  The union $\Xhat=\coprod_P X_P$
(with the quotient topology from the natural map $\Xbar\to \Xhat$) is the
{\itshape reductive Borel-Serre compactification\/}; it was introduced by
Zucker \cite{refnZuckerWarped}.  The reductive Borel-Serre
compactification is natural from a differential geometric standpoint since
the locally symmetric metric on $X$ degenerates precisely along these
nilmanifolds near the boundary of $\Xbar$.

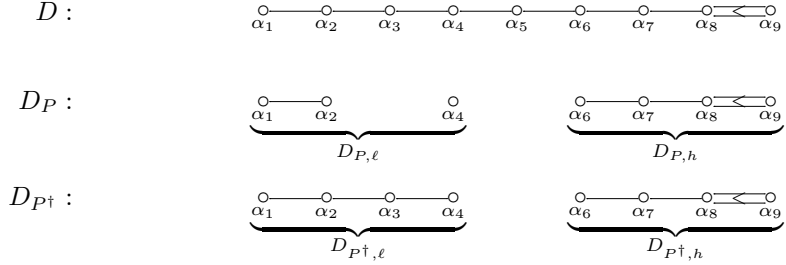
\begin{figure}
\entrymodifiers={}
\begin{align*}
&\xymatrix @!0 @M=0pt {{\txt{\llap{$D:$}}} &{}&{}&
{\circ}\save[]+<0in,-.075in>*{\scriptstyle\al_1}\restore \ar@{-}[r] &
{\circ}\save[]+<0in,-.075in>*{\scriptstyle\al_2}\restore \ar@{-}[r] &
{\circ}\save[]+<0in,-.075in>*{\scriptstyle\al_3}\restore \ar@{-}[r] &
{\circ}\save[]+<0in,-.075in>*{\scriptstyle\al_4}\restore \ar@{-}[r] &
{\circ}\save[]+<0in,-.075in>*{\scriptstyle\al_5}\restore \ar@{-}[r] &
{\circ}\save[]+<0in,-.075in>*{\scriptstyle\al_6}\restore \ar@{-}[r] &
{\circ}\save[]+<0in,-.075in>*{\scriptstyle\al_7}\restore \ar@{-}[r] &
{\circ}\save[]+<0in,-.075in>*{\scriptstyle\al_8}\restore
\ar@{-}@<.4ex>[r] \ar@{-}@<-.4ex>[r] \ar@{}[r]|{<}&
{\circ}\save[]+<0in,-.075in>*{\scriptstyle\al_9}\restore} \\[4ex]
&\xymatrix @!0 @M=0pt {{\txt{\llap{$D_P:$}}} &{}&{}&
{\circ}\save[]+<0in,-.075in>*{\scriptstyle\al_1}="a"\restore \ar@{-}[r] &
{\circ}\save[]+<0in,-.075in>*{\scriptstyle\al_2}\restore &
{} &
{\circ}\save[]+<0in,-.075in>*{\scriptstyle\al_4}="b"\restore
\save "a";"b" **{}?."a"."b" *!/u1ex/\frm{_\}} *!/u3.25ex/{\scriptstyle
D_{P,\l}}\restore &
{} &
{\circ}\save[]+<0in,-.075in>*{\scriptstyle\al_6}="a"\restore \ar@{-}[r] &
{\circ}\save[]+<0in,-.075in>*{\scriptstyle\al_7}\restore \ar@{-}[r] &
{\circ}\save[]+<0in,-.075in>*{\scriptstyle\al_8}\restore
\ar@{-}@<.4ex>[r] \ar@{-}@<-.4ex>[r] \ar@{}[r]|{<}&
{\circ}\save[]+<0in,-.075in>*{\scriptstyle\al_9}="b"\restore
\save "a";"b" **{}?."a"."b" *!/u1ex/\frm{_\}} *!/u3.25ex/{\scriptstyle
D_{P,h}}\restore} \\[1ex]
&\xymatrix @!0 @M=0pt {{\txt{\llap{$D_{P^\dag}:$}}} &{}&{}&
{\circ}\save[]+<0in,-.075in>*{\scriptstyle\al_1}="a"\restore \ar@{-}[r] &
{\circ}\save[]+<0in,-.075in>*{\scriptstyle\al_2}\restore \ar@{-}[r] &
{\circ}\save[]+<0in,-.075in>*{\scriptstyle\al_3}\restore \ar@{-}[r] &
{\circ}\save[]+<0in,-.075in>*{\scriptstyle\al_4}="b"\restore
\save "a";"b" **{}?."a"."b" *!/u1ex/\frm{_\}} *!/u3.25ex/{\scriptstyle
D_{P^\dag,\l}}\restore &
{} &
{\circ}\save[]+<0in,-.075in>*{\scriptstyle\al_6}="a"\restore \ar@{-}[r] &
{\circ}\save[]+<0in,-.075in>*{\scriptstyle\al_7}\restore \ar@{-}[r] &
{\circ}\save[]+<0in,-.075in>*{\scriptstyle\al_8}\restore
\ar@{-}@<.4ex>[r] \ar@{-}@<-.4ex>[r] \ar@{}[r]|{<}&
{\circ}\save[]+<0in,-.075in>*{\scriptstyle\al_9}="b"\restore
\save "a";"b" **{}?."a"."b" *!/u1ex/\frm{_\}} *!/u3.25ex/{\scriptstyle
D_{P^\dag,h}}\restore}
\end{align*}
\caption{An example of $D_P=D_{P,\l}\times D_{P,h}$ and
$D_{P^\dag}=D_{P^\dag,\l}\times D_{P^\dag,h}$}
\label{figFactorizationExample}
\end{figure}

Finally assume now that $D$ is Hermitian symmetric.  Then each $D_P$
factors into a product $D_{P,\l}\times D_{P,h}$, where $D_{P,h}$ is
again Hermitian symmetric (see Figure~\ref{figFactorizationExample}).
This induces a factorization (modulo a finite quotient) $X_P =
X_{P,\l}\times X_{P,h}$ of the arithmetic quotients and we consider
the projection $X_P\to X_{P,h}$ onto the second factor.  Now among the
different $P\in\Pl$ that yield the same $X_{P,h}$, let $P^\dag\in
\Pl_1$ be the maximal one and set $F_{P^\dag} = X_{P,h}$.  Thus each
stratum of $\Xhat$ has a projection $X_P\to F_{P^\dag}$.  The union
$\Xstar= \coprod_{R\in\Pl_1} F_R$ (with the quotient topology from the
map $\Xhat\to \Xstar$) is the {\itshape Baily-Borel-Satake
compactification\/} $\Xstar$.  Topologically $\Xstar$ was constructed
by Satake \cite{refnSatakeCompact}, \cite{refnSatakeQuotientCompact}
(though the description we have given is due to Zucker
\cite{refnZuckerSatakeCompactifications}); if $\G$ is contained in
the group of biholomorphisms of $D$, the compactification $\Xstar$ was
given the structure of a normal projective algebraic variety by Baily and
Borel \cite{refnBailyBorel}.

The simplest example where all three compactifications are distinct is
the Hilbert modular surface case.  Here $G=R_{k/\QQ} \SL(2)$ where $k$
is a real quadratic extension.  There is only one proper parabolic
$\QQ$-subgroup $P$ up to $G(\QQ)$-conjugacy; $Y_P$ is a torus bundle
over $X_P=S^1$ and $F_P$ is a point.

\section{The conjecture}
\label{sectRapoportsConjecture}
Assume that $D$ is Hermitian symmetric.  Let $E\in\Rep(G)$, the category of
finite dimensional regular representations of $G$ and let $\EE$ denote the
corresponding local system on $X$.  Let $\IC(\Xhat;\EE)$ and
$\IC(\Xstar;\EE)$ denote middle perversity intersection cohomology sheaves%
\footnote{By a ``sheaf'' we will always mean a complex of sheaves
representing an element of the derived category.  A derived functor
will be denoted by the same symbol as the original functor, thus we
will write $\pi_*$ instead of $R\pi_*$.}
on $\Xhat$ and $\Xstar$ respectively \cite{refnGoreskyMacPhersonIHII}.

For example, $\IC(\Xhat;\EE) = \t^{\leqslant p(\codim X_P)} j_{P*} \EE$ if
$\Xhat$ has only one singular stratum $X_P$; here $j_{P*}$ denotes the
derived direct image functor of the inclusion $j_P:\Xhat\setminus
X_P\hookrightarrow \Xhat$, $\codim X_P$ denotes the topological
codimension, $p(k)$ is one of the middle perversities
$\left\lfloor(k-1)/2\right\rfloor$ or $\left\lfloor(k-2)/2\right\rfloor$,
and $\t^{\leqslant p(k)}$ truncates link cohomology in degrees $> p(k)$.
In general the pattern of pushforward/truncate is repeated over each
singular stratum.  Note that since $\Xhat$ may have odd codimension strata,
$\IC(\Xhat;\EE)$ depends on the choice of the middle perversity $p$; on the
other hand, since $\Xstar$ only has even codimension strata,
$\IC(\Xstar;\EE)$ is independent of $p$.

\begin{maintheorem}
Let $X$ be  an arithmetic
quotient of a Hermitian symmetric space.  Then $\pi_* \IC(\Xhat;\EE)
\cong \IC(\Xstar;\EE)$.  \textup(That is, they are isomorphic in the derived
category.\textup)
\end{maintheorem}

Following discussions with Kottwitz, Rapoport conjectured the theorem in a
letter to Borel \cite{refnRapoportLetterBorel} and later provided
motivation for it in an unpublished note \cite{refnRapoportNote}.
Previously Zucker had noticed that the conjecture held for $G=\Sympl(4)$,
$\EE=\CC$.  The conjecture was later rediscovered by Goresky and MacPherson
and described in an unpublished preprint
\cite{refnGoreskyMacPhersonWeighted} in which they also announced the
theorem for $G=\Sympl(4)$, $\Sympl(6)$, and (for $\EE=\CC$) $\Sympl(8)$.
The first published appearance of the conjecture was in a revised version
of Rapoport's note \cite{refnRapoport} and included an appendix by Saper
and Stern giving a proof of the theorem when $\QQrank G=1$.

To see one reason why the conjecture might be useful in the theory of
automorphic forms, note that the right hand side $\IC(\Xstar;\EE)$ is
isomorphic to the $L^2$-cohomology sheaf $\L_{(2)}(\Xstar;\EE)$ by (the proof
of) Zucker's conjecture \cite{refnLooijenga},
\cite{refnSaperSternTwo}.  The trace of a Hecke
operator on $L^2$-cohomology could then be studied topologically via the
Lefschetz fixed point formula for $\IC(\Xstar;\EE)$.  However the
singularities of $\Xhat$ are simpler than those of $\Xstar$ so a Lefschetz
fixed point formula for $\IC(\Xhat;\EE)$ should be easier to calculate.
The conjecture says that this should give the same result.  Also note that
a Lefschetz fixed point formula for $\IC(\Xhat;\EE)$ involves a sum over
$\Pl$, while a Lefschetz fixed point formula for $\IC(\Xstar;\EE)$ involves
a sum over $\Pl_1$.  Thus it is more likely that the former can be directly
related to the Arthur-Selberg trace formula for a Hecke operator on
$L^2$-cohomology \cite{refnArthur}.

This program has been pursued by Goresky and MacPherson, but instead of
$\IC(\Xhat;\EE)$ they use the ``middle weighted cohomology''
$\WC(\Xhat;\EE)$ in which cohomology classes in the link are truncated
according to their {\itshape weight\/} as opposed to their degree.  Thus
weighted cohomology is an algebraic analogue of $L^2$-cohomology.  Goresky
and MacPherson prove (in joint work with Harder
\cite{refnGoreskyHarderMacPherson}) the analogue of the above theorem,
$\pi_* \WC(\Xhat;\EE) \cong \IC(\Xstar;\EE)$, calculate the Lefschetz fixed
point formula \cite{refnGoreskyMacPhersonTopologicalTraceFormula}, and (in
joint work with Kottwitz) show that it agrees with Arthur's trace formula
for $L^2$-cohomology \cite{refnGoreskyKottwitzMacPherson}.

Nonetheless the original conjecture remains interesting for a number of
reasons.  First of all, intersection cohomology is a true topological
invariant and the local cohomology of $\IC(\Xhat;\EE)$ behaves better than
that of $\WC(\Xhat;\EE)$ when $\EE$ varies.  Secondly, the local property
(``micro-purity'') one needs to prove is much deeper for $\IC(\Xhat;\EE)$
than for $\WC(\Xhat;\EE)$ and should have applications elsewhere.  And
finally the method used to attack the conjecture, the theory of
$\L$-modules, has application to other cohomology, in particular, weighted
cohomology, $L^2$-cohomology, and ordinary cohomology.

In \S\S\ref{sectLModules}--\ref{sectMicroSupportFunctoriality}\ we will
indicate how the Main Theorem follows from three theorems in the
theory of $\L$-modules.

\section{A generalization}
\label{sectGeneralization}
This section is optional; we will indicate a more general context in which
the Main Theorem holds.  First we sketchily recall the general theory of
Satake compactifications \cite{refnSatakeCompact},
\cite{refnSatakeQuotientCompact}, \cite{refnZuckerSatakeCompactifications},
\cite{refnCasselmanGeometricRationality}.  By embedding $D$ into a real
projective space via a finite-dimensional representation $\sigma$ of $G$
and then taking the closure, Satake constructed a finite family of {\it
Satake compactifications\/} $\lsb\RR \Dstar$ of $D$.  Each of these is
equipped with an action of $G(\RR)$ and is formed by adjoining to $D$
certain {\itshape real boundary components\/}.  Let $\Dstar$ denote the
union of $D$ together with those real boundary components whose normalizer
is defined over $\QQ$; call these the {\itshape rational boundary
components\/}.  In the {\itshape geometrically rational\/} case (a
condition satisfied for example if $\sigma$ is $\QQ$-rational%
\footnote{Borel points out that in his 1962 Bruxelles conference paper
  ``Ensembles fondamentaux pour les groupes arithm{\'e}tiques'' he proves
  geometric rationality only when $\sigma$ is strongly $\QQ$-rational.  In
  \cite{refnSaperGeometricRationality} we prove geometric rationality for
  the general $\QQ$-rational case.}%
) one may
equip $\Dstar$ with a suitable topology so that $\Xstar=\G\back \Dstar$ is
a Hausdorff compactification of $X$; this is also called a {\itshape Satake
compactification\/}.  For $D$ Hermitian symmetric, one of the Satake
compactifications is (topologically equivalent to) the closure of the
realization of $D$ as a bounded symmetric domain and it is geometrically
rational; the corresponding compactification of $X$ is the
Baily-Borel-Satake compactification.

Let $\lsp0 G = \bigcap_{\chi\in X_\QQ(G)} \Ker \chi^2$ so that $G(\RR) =
\lsp0 G(\RR)A_G$ \cite{refnBorelSerre}.  Suppose that $\rank \lsp0 G =
\rank K$, that is, $\lsp0 G(\RR)$ has discrete series representations.
This is equivalent to the assumption that the maximal $\RR$-split torus in
the center of $G$ is also $\QQ$-split and that the real points of
$G^{\text{der}}$ (the semisimple derived group) has discrete series
representations.  (We may also substitute here the adjoint group
$G^{\text{ad}}$ for $G^{\text{der}}$.)  We say in this case that $D$ is an
{\itshape equal-rank symmetric space}.  A Satake compactification $\lsb\RR
\Dstar$ of $D$ will be called a {\itshape real equal-rank Satake
compactification\/} if all the real boundary components of $\lsb\RR \Dstar$
are also equal-rank symmetric spaces.  The possible $D$ that admit real
equal-rank Satake compactifications are listed in
\cite{refnZuckerLtwoIHTwo}; they include the Hermitian symmetric cases but
there are other infinite families as well.  If such a $\lsb\RR \Dstar$ is
geometrically rational%
\footnote{We show in \cite{refnSaperGeometricRationality} that this always
  holds except for certain explicitly described situations in $\QQrank$ $1$
  and $2$ involving restriction of scalars.}
then the corresponding compactification $\Xstar$ of
$X$ is also called a {\itshape real equal-rank Satake compactification\/};
note that we impose the equal-rank condition on all real boundary
components even though only the rational boundary components contribute to
$\Xstar$.

The generalization we alluded to above is that the Main Theorem holds for
real equal-rank Satake compactifications.  (Note that Borel conjectured
that the analogue of the Zucker conjecture should remain true for such
$\Xstar$ and Saper and Stern (unpublished) observed that their proof could be
adapted to this case.)

\section{$\L$-modules}
\label{sectLModules}
Now again let $G$ be any connected reductive group over $\QQ$ (with no
Hermitian hypothesis).  The ``sheaf'' $\IC(\Xhat;\EE)$ is actually an
object of $\Der_\X(\Xhat)$, the derived category of complexes of sheaves
$\Sheaf$ on $\Xhat$ that are constructible.  Here the constructibility of
$\Sheaf$ means that if for all $P\in \Pl$ we let $i_P:X_P\hookrightarrow
\Xhat$ denote the inclusion, then the local cohomology sheaf
$H(i_P^*\Sheaf)=H(\Sheaf|_{X_P})$ is locally constant, or equivalently the
cohomology sheaf $\EE_P = H(i_P^!\Sheaf)$ is locally constant on $X_P$.
Thus by the correspondence between local systems and representations of the
fundamental group one obtains a family of objects $E_P\in
\Graded(\G_{L_P})$, the category of graded $\G_{L_P}$-modules, one for each
$P\in\Pl$.

Instead of $\Sheaf$ we wish to work with a combinatorial analogue in
which $\Graded(\G_{L_P})$ is replaced by $\Graded(L_P)$, the category of
graded regular $L_P$-modules.  This analogue is what we will call an {\it
$\L$-module on $\Xhat$\/}.  We will describe just what an $\L$-module is
more precisely later, but first let us give some of the properties of the
categories $\Rep(\L_W)$ of $\L$-modules on $W$, where $W$ is any locally
closed union of strata of $\Xhat$:
\begin{enumerate}
\item if $W=X_P$, then $\Rep(\L_{X_P})= \Complex(L_P)$, the category of
complexes of regular $L_P$-modules;
\label{LModulesPropertiesOne}
\item for any inclusion $j:W\hookrightarrow W'$, there exist functors
$j^*,j^!: \Rep(\L_{W'})\to \Rep(\L_W)$ and $j_*,j_!:\Rep(\L_W)\to
\Rep(\L_{W'})$, as well as a degree truncation functor $\t^{\leqslant p}:
\Rep(\L_W)\to \Rep(\L_W)$;
\label{LModulesPropertiesTwo}
\item there is a realization functor $\Sheaf_W:\Rep(\L_W)\to \Der_\X(W)$ which
commutes with the functors in (\ref{LModulesPropertiesTwo}) and for
which the following diagram commutes:%
\entrymodifiers={++!}%
\newbox\fakeone%
\setbox\fakeone=\hbox{$\Rep(\L_{X_P})$}%
\newbox\faketwo%
\setbox\faketwo=\hbox{$\Graded(\G_{L_P})$}%
\begin{equation*}
\xymatrix@C+=5ex{
{\Rep(\L_{X_P})}  \ar@<1ex>[d]^{H(\cdot)} \ar@<.5ex>[r] & *++!\hbox to
\wd\faketwo{\hfill$\Der_\X(X_P)$\hfill}
\ar@<1ex>[d]^{H(\cdot)} \\
*++!\hbox to \wd\fakeone{\hfill$\Graded(L_P)$\hfill} \ar@<1ex>[u]
\ar@<.5ex>[r]^-{\Res} &
{\Graded(\G_{L_P})\rlap{\quad .}}
\ar@<1ex>[u]
}
\end{equation*}
\label{LModulesPropertiesThree}
\end{enumerate}
Note that one advantage of $\L$-modules over sheaves is that the left hand
vertical arrows in (\ref{LModulesPropertiesThree}) are equivalences of
categories, unlike those on the 
right; this is because $\Rep(L_P)$ is a semisimple category.

So roughly speaking an $\L$-module is like a sheaf $\Sheaf$ with the ``extra
structure'' that $\EE_P = H(i_P^!\Sheaf)$ is associated to a regular
$L_P$-module, as opposed to merely a $\G_{L_P}$-module.  Condition
(\ref{LModulesPropertiesTwo})
implies that the usual operations on sheaves preserve this ``extra
structure''.  The following example shows this is reasonable.  Let $\EE$ be
a local system on $X$ associated to a regular representation $E$ of $G$.
The smooth part of the link bundle of a real codimension $k$ stratum
$X_P\subset \Xhat$ is the flat bundle with fiber $|\D_P|^\circ\times
\G_{N_P}\back N_P(\RR)$, where $|\D_P|^\circ$ is an open $(k-1)$-simplex and
$\G_{L_P}$ acts via conjugation on the second factor
\cite[\S8]{refnGoreskyHarderMacPherson}.  Thus $H(i_P^*i_{G*}\EE)
\cong \HH(\G_{N_P}\back N_P(\RR);\EE)$, the local system associated to the
$\G_{L_P}$-module $H(\G_{N_P}\back N_P(\RR);\EE)$.  However by van Est's
theorem \cite{refnvanEst}, $H(\G_{N_P}\back N_P(\RR);\EE)$ is isomorphic to
the restriction of the regular $L_P$-module $H(\n_P;E)$, where $\n_P$ is
the Lie algebra of $N_P(\RR)$.

In fact this also suggests how to precisely define $\L$-modules.  Let
$\Pl(W)\subseteq \Pl$ correspond to the strata of $W$.  For $P\le Q$ let
$\n_P^Q$ be the Lie algebra of $N_P(\RR)/N_Q(\RR)$.  An {\itshape $\L$-module\/}
$\M\in \Rep(\L_W)$ is a family $(E_\cdot,f_{\cdot\cdot})$ consisting of
objects $E_P\in\Graded(L_P)$ for every $P\in\Pl(W)$ and degree 1 morphisms
$f_{PQ}:H(\n_P^Q;E_Q)\xrightarrow{[1]} E_P$ for every $P\le Q\in
\Pl(W)$ such that
\begin{equation*}
\sum_{P\le Q\le R} f_{PQ}\circ H(\n_P^Q;f_{QR}) = 0
\end{equation*}
for all $P\le R\in\Pl(W)$.  The functors $i_P^!$ and $i_P^*$ are given
by
\begin{align*}
i_P^!\M &= (E_P,f_{PP})\ , \\
i_P^*\M &= \biggl( \bigoplus_{P\le R}
H(\n_P^R;E_R), \sum_{P\le R\le S} H(\n_P^R;f_{RS}) \biggr)\ .
\end{align*}

We define the global cohomology $H(\Xhat;\M)$ of an $\L$-module $\M$ to be
the hypercohomology of its realization, $H(\Xhat;\Sheaf_{\Xhat}(\M))$.  In
general we will often write simply $\M$ for both the $\L$-module and its
realization $\Sheaf_{\Xhat}(\M)$; it should be clear what is meant from the
context.

\section{Examples of $\L$-modules}
\label{sectExamples}
\begin{enumerate}
\item Let $E\in \Rep(G)$.  Then the $\L$-module $i_{G*}E$ defined by
$E_G=E$ and $E_P=0$ for $P\neq G$ corresponds via $\Sheaf_{\Xhat}$ to
$i_{G*}\EE$ and its cohomology is the ordinary cohomology
$H(X;\EE)=H(\G;E)$.
\label{ExampleOne}

\item It follows immediately from the properties of $\L$-modules in the
previous section that given $E\in \Rep(G)$ there exists an $\L$-module
$\IC(\Xhat;E)$ which maps under $\Sheaf_{\Xhat}$ to the intersection
cohomology sheaf $\IC(\Xhat;\EE)$.  For example, if $\Pl=\{G,P\}$ (that is,
$\Xhat$ has only one singular stratum) and $p=p(\codim X_P)$, then
\begin{equation*}
\IC(\Xhat;E)= 
\begin{pmatrix} E_G=E, \  E_P = (\t^{>p}H(\n_P;E))[-1], \\
 f_{PG} : H(\n_P;E) \to \t^{>p} H(\n_P;E)
\end{pmatrix}
\end{equation*}
where $\t^{>p} H(\n_P;E) = \bigoplus_{i>p} H^i(\n_P;E)[-i]$ and
$f_{PG}$ is the projection.  Note that the truncation $\t^{\leqslant
p}$ of local cohomology at $X_P$ has been implemented {\it
externally\/} via a mapping cone; this is valid in view of the
quasi-isomorphism $\t^{\leqslant p} C \tildearrow \Cone(C\!\to\! \t^{>
p}C)[-1]$ for any complex $C$.

\item The weighted cohomology sheaf and the $L^2$-cohomology sheaf may also
be lifted to $\L$-modules $\WC(\Xhat;E)$ and $\L_{(2)}(\Xhat;E)$; for the
latter we must replace $\Rep(L_P)$ by the category of {\itshape locally
regular\/} $L_P$-modules to handle the potentially infinite dimensional
local cohomologies.
\end{enumerate}

\section{Micro-support of $\L$-modules}
\label{sectMicroSupport}
The {\itshape support\/} of a sheaf $\Sheaf$ is the set of points $x$ such that
$H(\Sheaf)_x\neq 0$.  As is well-known the global cohomology of $\Sheaf$
vanishes if the support is empty (that is, the sheaf is
quasi-isomorphic to $0$).  For an $\L$-module $\M$ we will state in
the next section a more subtle vanishing result based on the {\it
micro-support\/} of $\M$ which we now define; this is a rough analogue
of the corresponding notion for sheaves \cite{refnKashiwaraSchapira}.

Let $P\in \Pl$ and let $\IrrRep(L_P)$ denote the set of irreducible regular
$L_P$-modules.  For $V\in \IrrRep(L_P)$ let $\xi_V$ be the character by
which $A_P$ acts on $V$.  Let $\D_P$ be the simple roots of the adjoint
action of $A_P$ on $\n_P$; the parabolic $\QQ$-subgroups $Q\ge P$ are
indexed by subsets $\D_P^Q$ of $\D_P$.  Define $P\le Q_V\le Q_V'\in \Pl$ by
\begin{align*}
\D_P^{Q_V} &=\{\,\al\in \D_P\mid (\xi_V+\r,\al)<0\,\} \ , \\
\D_P^{Q_V'} &=\{\,\al\in \D_P\mid (\xi_V+\r,\al)\le 0\,\}\ ,
\end{align*}
where $\r$ denotes one-half the sum of the positive roots of $G$ and the
inner product is induced by the Killing form of $G$.  Let $M_P = \lsp0 L_P$
so that $L_P(\RR) = M_P(\RR)A_P$.  Let $V|_{M_P}$ denote the restriction of
the representation $V$  to $M_P$.

The {\itshape micro-support $\mS(\M)$ of $\M$} is the subset of
$\coprod_{P\in\Pl} \IrrRep(L_P)$ consisting of those $V\in \IrrRep(L_P)$
satisfying
\begin{enumerate}
\item $(V|_{M_P})^* \cong \overline{V|_{M_P}}$, and
\label{SScondOne}
\item there exists $Q_V\le Q \le Q_V'$ such that
\label{SScondTwo}
\begin{equation}
H(i_P^* \ihat_Q^! \M)_V \neq 0\ .
\label{eqnGeneralType}
\end{equation}
\end{enumerate}
Here $\ihat_Q: \Xhat_Q \hookrightarrow \Xhat$ is the inclusion of the
closure of the stratum $X_Q$ and the subscript $V$ indicates the
$V$-isotypical component.  A simple example of the computation of
micro-support will be given in \S\ref{sectTilouine}.

Condition ~(\ref{SScondOne}) is equivalent to the existence of a
nondegenerate sesquilinear form on $V$ which is invariant under the action
of $M_P$.

As for condition (\ref{SScondTwo}), let $\jhat_Q:\Xhat \setminus
\Xhat_Q \hookrightarrow \Xhat$ be the open inclusion.  Note that we have a short exact sequence
\begin{equation*}
0\to i_P^* \ihat_Q^! \M \to i_P^* \M  \to i_P^* \jhat_{Q*}\jhat_Q^*\M \to 0
\end{equation*}
and a corresponding long exact sequence.  Topologically, this is the
long exact sequence of the pair $(U,U\setminus(U\cap \Xhat_Q))$ where
$U$ is a small neighborhood of a point of $X_P$.  Thus condition (\ref{SScondTwo})
means that
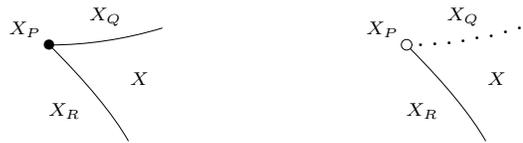
\begin{equation*}
\xy\xycompile{ 0;<.75mm,0mm>:
(0,20)*+!DR{\scriptstyle X_P};(20,23) **\crv{(10,20)},
(10,22)*+!D{\scriptstyle X_Q},
(0,20)*-{\bullet};(14,3) **\crv{(10,10)},
(7,11)*+!RU{\scriptstyle X_R},(16,14)*+{\scriptstyle X},
(10,-5)*{H(U;\M)_V}
}\endxy\qquad\xy 0;<.75mm,0mm>:
(5,-5)\ar @{>}(15,-5)
\endxy\qquad\xy\xycompile{ 0;<.75mm,0mm>:
(0,20)*+!DR{\scriptstyle X_P};(20,23) **\crv{~*=<4pt>{.} (10,20)},
(10,22)*+!D{\scriptstyle X_Q},
(0,20)*\cir<2pt>{};(14,3) **\crv{(10,10)},
(7,11)*+!RU{\scriptstyle X_R},(16,14)*+{\scriptstyle X},
(10,-5)*{H(U\setminus(U\cap \Xhat_Q);\M)_V}}
\endxy
\end{equation*}
is not an isomorphism for some degree and for some $Q$ between $Q_V$ and
$Q_V'$.

It is convenient to define the {\itshape essential micro-support
$\emS(\M)$ of $\M$\/} to be the subset consisting of those $V\in
\mS(\M)$ for which
\begin{equation*}
\Type_V(\M) = \Im \bigl( H(i_P^* \ihat_{Q_V}^! \M)_V  \longrightarrow H(i_P^*
\ihat_{Q_V'}^! \M)_V \bigr)
\end{equation*}
is nonzero.  The essential micro-support of $\M$ determines the
micro-support (though not the actual parabolics $Q$ that arise in
condition (\ref{SScondTwo})).  In fact the relation between $\mS(\M)$ and $\emS(\M)$
is analogous to the relation between the strata of a nonreduced
variety (possibly with embedded components) and the smooth open strata
of the irreducible components: there exists a partial order $\preccurlyeq$ on
$\coprod_{P\in\Pl} \IrrRep(L_P)$ such that if $V\in \mS(\M)$ then
there exists $\tilde V\in \emS(\M)$ with $V\preccurlyeq \tilde V$, and if
$\tilde V\in \emS(\M)$ and $V\preccurlyeq \tilde V$ then $V\in\mS(\M)$.

\section{A vanishing theorem for $\L$-modules}
\label{sectVanishingTheorem}
The justification for the definition of $\mS(\M)$ is that it is an
ingredient for a vanishing theorem for $H(\Xhat;\M)$.  To state the
theorem we need some more notation.

Let $V\in\IrrRep(L_P)$ have highest weight $\u\in \h_\CC^*$ where $\h$ is a
fundamental (maximally compact) Cartan subalgebra for the Lie algebra
$\levi_P$ of $L_P(\RR)$ equipped with a compatible ordering.  Assume
$(V|_{M_P})^* \cong \overline{ V|_{M_P}}$ and define
\begin{align*}
L_P(\u) &= \text{ the centralizer of $\u\in\h_\CC^*\subset
\levi_{P\CC}^*$\ ,} \\
&= \text{ the reductive subgroup of $L_P$ with
roots $\{\,\g\in\Phi(\levi_{P\CC},\h_\CC)\mid (\g,\u)=0\,\}$\ ,} \\
D_P(\u) &= \text{ the associated symmetric space } L_P(\u)(\RR)/(K_P\cap
L_P(\u))A_P\ .
\end{align*}
Choose a compatible ordering for which $\dim D_P(\u)$ is maximized and let
$D_P(V) = D_P(\u)$.  Suppose now that $V\in\emS(\M)$.  Let $c(V;\M)\le
d(V;\M)$ be the least and greatest degrees in which $\Type_V(\M)$ is
nonzero, and define
\begin{align*}
\tilde c(V;\M) &=   \tfrac12(\dim D_P - \dim D_P(V)) + c(V;\M) \ , \\
\tilde d(V;\M) &=   \tfrac12(\dim D_P + \dim D_P(V)) + d(V;\M) \ .
\end{align*}
Set
\begin{equation*}
c(\M) = \inf_{V\in\emS(\M)} \tilde c(V;\M)\ ,\qquad d(\M) =
\sup_{V\in\emS(\M)} \tilde d(V;\M)\ .
\end{equation*}
(One can show that the same values are obtained if instead we consider all
$V\in\mS(\M)$ and let $c(V;\M)\le d(V;\M)$ be the least and greatest
degrees in which \eqref{eqnGeneralType} is nonzero (for any $Q$).)

\begin{thm}
\label{thmVanishing}
$H^i(\Xhat;\M) = 0$ for $i\notin [c(\M),d(\M)]$.
\end{thm}

Let us comment briefly on the proof which uses combinatorial Hodge-de Rham
theory.  The sheaf $\Sheaf_{\Xhat}(\M)$ has an incarnation as a complex of
fine sheaves whose global sections are ``combinatorial'' differential
forms.  That is, an element of $\G(\Xhat;\Sheaf_{\Xhat}^\cdot(\M))$ is a
family $(\o_P)_{P\in\Pl}$, where each $\o_P$ is a special differential form
on $X_P$ with coefficients in $\EE_P$.  (For $P=G$, the {\itshape special
differential forms\/} \cite[(13.2)]{refnGoreskyHarderMacPherson} on $X=X_G$
are those which near each boundary stratum $Y_Q$ of the Borel-Serre
compactification $\Xbar$ are the pullback of an $N_Q(\RR)$-invariant form
on $Y_Q$; they form a resolution of $\EE_G$.)  The differential is a sum of
the usual de Rham exterior derivative (on each $\o_P$) together with
operators based on the $f_{PQ}$.

To do harmonic theory we need a metric; unfortunately the locally
symmetric metric on each $X_P$ is not appropriate since it would
introduce unwanted $L^2$-growth conditions on the differential forms.
Instead the theory of {\itshape tilings\/} from \cite{refnSaperTilings}
gives a natural piecewise analytic diffeomorphism of $\Xbar$ onto a
closed subdomain $\Xbar_0$ of the interior $X$; the pullback of the
locally symmetric metric under this map yields metrics on all $X_P$
which extend to nondegenerate metrics on their boundary strata.  Now a
spectral analogue of the Mayer-Vietoris sequence as in
\cite{refnSaperSternTwo} reduces the problem to a vanishing theorem
for combinatorial $L^2$-cohomology near each stratum $X_P$.  After
unraveling the combinatorics one obtains contributions to the
cohomology of the form $H_{(2)}(X_P;\VV)\otimes \Type_V(\M)$ for $V
\in \emS(\M)$; by Raghunathan's vanishing theorem
\cite{refnRaghunathan}, \cite{refnRaghunathanCorrection},
\cite{refnSaperSternTwo} this is zero outside the degree range
$[\tilde c(V;\M),\tilde d(V;\M)]$.  (The proof is actually more complicated since
there are infinite dimensional contributions from $\mS(\M)\setminus
\emS(\M)$ as well.)

\section{Micro-purity of intersection cohomology}
\label{sectIHMicroPurity}
We will say an $\L$-module $\M$ on $\Xhat$ is {\itshape $V$-micro-pure\/} if
$\emS(\M)=\{V\}$ with $\Type_V(\M)$ concentrated in degree $0$.

\begin{thm}
\label{thmMicroPurity}
Assume the irreducible components of the $\QQ$-root
system of $G$ are of type $A_n$, $B_n$, $C_n$, $BC_n$, or $G_2$.  Let
$E\in\IrrRep(G)$ satisfy $(E|_{\lsp0 G})^* \cong \overline{ E|_{\lsp0 G}}$.
Then $\IC(\Xhat;E)$ is $E$-micro-pure.
\end{thm}

If $D$ is a Hermitian symmetric space (or an equal-rank symmetric space
admitting a real equal-rank Satake compactification as in
\S\ref{sectGeneralization}) $G$ will have a $\QQ$-root system of the indicated
type and thus the theorem applies in the context of Rapoport and
Goresky-MacPherson's conjecture.  In fact it is quite possible that this
restriction in the theorem may be removed; it is only required at one
crucial stage in the proof.

What the theorem is asserting is that $V\notin \emS(\IC(\Xhat;E))$ for
$V\in\IrrRep(L_P)$ with $P\neq G$.  When $P$ is a maximal parabolic we
can give a brief indication of how this is proven; for definiteness we
assume $p$ is the upper middle perversity.  In this case
\begin{equation}
H(i_P^*\ihat_Q^!\IC(\Xhat;E)) =
\begin{cases} \t^{\leqslant p} H(\n_P;E) & \text{for $Q=G$\ ,} \\
              (\t^{>p} H(\n_P;E))[-1]	 & \text{for $Q=P$\ ,}
\end{cases}
\label{eqnMicroTypes}
\end{equation}
where $p= \lfloor\tfrac12 \dim\n_P\rfloor$.  Let $\lambda$ be the
highest weight of $E$.  By Kostant's theorem \cite{refnKostant} an
irreducible component $V$ of $H(\n;E)$ has highest weight
$w(\lambda+\r)-\r$ where
\begin{equation*}
w \in W_P=\{w\in W\mid w^{-1}\g>0 \text{ for all postive roots $\g$ of
$\levi_{P\CC}$}\}\ ,
\end{equation*}
the set of minimal length representatives of the Weyl group quotient
$W_{L_P}\back W$.  Furthermore $V$ occurs in degree $\l(w)$, the
length of $w$, with multiplicity $1$.  Assume now that $V\in
\emS(\IC(\Xhat;E))$.  Since the two cases in \eqref{eqnMicroTypes}
above do not share a common component we must have $Q_V=Q_V'$, that
is, $(\xi_V+\r,\al)\neq 0$ for the unique $\al\in \D_P$.  Furthermore
\eqref{eqnMicroTypes} also shows that the possibilities
$(\xi_V+\r,\al) < 0$ and $(\xi_V+\r,\al) > 0$ correspond respectively
to $\l(w) \le \tfrac12\dim \n_P$ and $\l(w) > \tfrac12\dim \n_P$.
However the following lemma from \cite{refnSaperLModules} shows that
in fact the opposite relation between weight and degree holds (the
nonnegative term $\dim \n_P(V)$ here may be ignored for now---it will be
defined in \S\ref{sectTilouine}):

\begin{lem}
\label{lemBasicLemma}
Let $V\in \IrrRep(L_P)$ have highest weight $w(\lambda+\r)-\r$ where
$w\in W_P$ and $\lambda\in\h_\CC^*$ is dominant.  Assume that
$(V|_{M_P})^* \cong \overline{ V|_{M_P}}$.
\begin{enumerate}
\item If $(\xi_V+\r,\al)\le 0$ for all $\al\in\D_P$, then
$\l(w) \ge \tfrac12(\dim \n_P + \dim \n_P(V))$.
\label{BasicLemmaPartOne}
\item If $(\xi_V+\r,\al)\ge 0$ for all $\al\in\D_P$, then
$\l(w) \le \tfrac12(\dim \n_P - \dim \n_P(V))$.
\end{enumerate}
\end{lem}

The only remaining possibility is that $\l(w) = \tfrac12\dim \n_P$, but
since $(\xi_V+\r,\al) \neq 0$ and $(E|_{\lsp0 G})^* \cong \overline{
E|_{\lsp0 G}}$ this is impossible by an argument based on
\cite{refnBorelCasselman}.  By the way, Lemma~\ref{lemBasicLemma} is basic
to the proofs of Theorems ~\ref{thmVanishing}, ~\ref{thmFunctoriality}, and
\ref{thmTilouine} as well and has its origin in a result of Casselman for
$\RRrank$ one \cite{refnCasselman}.

When $P$ is not a maximal parabolic the situation is far more complicated.
The irreducible components of $H(i_P^*\IC(\Xhat;E))$ are among those of
$H(i_P^*i_{G*}i_G^*\IC(\Xhat;E)) = H(\n_P;E)$, but they may occur in
various degrees and with multiplicity.  Since we do not know a nonrecursive
formula for $H(i_P^*\IC(\Xhat;E))$ we must rely on the inductive
definition.  However condition (\ref{SScondOne}) in the definition of
micro-support is not preserved upon passing to a larger stratum.
Specifically, let $P< R$ and suppose $V$ is an irreducible component of
$H(\n_P;E)= H(\n_P^R; H(\n_R;E))$.  It must lie within $H(\n_P^R;V_R)$ for
some irreducible component $V_R$ of $H(\n_R;E)$.  The difficulty in using
induction is that $(V|_{M_P})^* \cong \overline{V|_{M_P}}$ does not
imply $(V_R|_{M_R})^* \cong \overline{V_R|_{M_R}}$.

These difficulties do not apply to $\WC(\Xhat;E)$ and in fact a fairly
simple argument shows that Theorem~\ref{thmMicroPurity} holds for
$\WC(\Xhat;E)$ without any hypothesis on the $\QQ$-root system and for
either middle weight profile.  Indeed since $\WC(\Xhat;E)$ is defined
directly in terms of weight the relationship between weight and degree
provided by Lemma~\ref{lemBasicLemma} is not needed and hence the condition
$(V|_{M_P})^* \cong \overline{V|_{M_P}}$ plays no role in the
proof.

\section{Functoriality of micro-support and proof of the Main Theorem}
\label{sectMicroSupportFunctoriality}
Let $\M$ be an $\L$-module which is $E$-micro-pure (for example,
$\M=\IC(\Xhat;E)$ by Theorem~\ref{thmMicroPurity}) and assume we are
in the context of
Rapoport and Goresky-MacPherson's conjecture, that is, $D$ is Hermitian
symmetric and $\pi:\Xhat\to \Xstar$ is the projection onto the
Baily-Borel-Satake compactification.  The desired equality $\pi_*\M=
\IC(\Xstar;\EE)$ is equivalent to certain local vanishing and covanishing
conditions on $\pi_*\M$ \cite{refnGoreskyMacPhersonIHII}.  To state them,
let $i_x:\{x\}\hookrightarrow \Xstar$ denote the inclusion of a point in a
stratum $F_R\subset \Xstar$.  Since every stratum of $\Xstar$ has even
codimension, $p(\codim F_R)= \tfrac12 \codim F_R -1$.  The local conditions
that characterize intersection cohomology now can be expressed as
\begin{equation}
\begin{aligned}
H^i(i_x^*\pi_*\M) &=0 \qquad \text{for $x\in F_R$, $i\ge\tfrac12\codim
F_R$, and }\\
H^i(i_x^!\pi_*\M) &=0 \qquad \text{for $x\in F_R$, $i\le\tfrac12\codim F_R$}
\end{aligned}
\label{VanishingCovanishingCondition}
\end{equation}
\begin{wrapfigure}{r}{25mm}
\begin{center}
\begin{xy}
<0mm,-5mm>;<1mm,-5mm>:
<5mm,-1mm>="c",
(15,0)="adown";(0,10)="bdown" **\crv{(10,5)&(5,8)}
?(.5)="xdown"*-{\scriptstyle \bullet},
"adown";+"c" **\crv{(17.5,0)},
"bdown";{"bdown"+"c"+<.4mm,.8mm>} **\crv{~*=<5pt>{.} "bdown"+<2.5mm,0mm>},
"adown"+"c";{"bdown"+"c"+<.4mm,.8mm>} **\crv{~*\dir{} (10,5)+"c" & (5,8)+"c"}
?(.33)="e";"adown"+"c" **@{-},
"e";{"bdown"+"c"+<.4mm,.8mm>} **\crv{~*=<5pt>{.} (10,5)+"c" & (5,8)+"c"},
"e"*++!LD{\scriptstyle X},
<0mm,5mm>;<1mm,5mm>:
<4mm,4mm>="d",
(15,0)="aup";(0,10)="bup" **\crv{(10,5)&(5,8)}
?(.5)="xup"*-{\scriptstyle \bullet},
"aup";+"d" **\crv{(17,0)},
"bup";{"bup"+"d"} **\crv{"bup"+<2mm,0mm>},
"aup"+"d";{"bup"+"d"} **\crv{"d"+(10,5) & "d"+(5,8)},
"adown";"aup" **@{-},
"bdown";"bup" **@{-},
"xdown";"xup" **@{-} ?(.5)*!RD{\scriptstyle \Xhat_{R,\l}},
<0mm,-25mm>;<1mm,-25mm>:
<5mm,-1mm>="c",<4mm,4mm>="d",
(15,0)="a";(0,10)="b"
**\crv{(10,5)&(5,8)}
?(.5)="x"*-{\scriptstyle \bullet}*+!UR{\scriptstyle x}
?(.2)*!UR{\scriptstyle F_R},
"a";+"c" **\crv{(17.5,0)},
"a";+"d" **\crv{(17,0)},
"b";{"b"+"c"+<.4mm,.8mm>} **\crv{~*=<5pt>{.} "b"+<2.5mm,0mm>},
"b";{"b"+"d"} **\crv{"b"+<2mm,0mm>},
"a"+"d";{"b"+"d"} **\crv{"d"+(10,5) & "d"+(5,8)},
"a"+"c";{"b"+"c"+<.4mm,.8mm>} **\crv{~*\dir{} "c"+(10,5) & "c"+(5,8)}
?(.17)="e";"a"+"c" **@{-},
"e";{"b"+"c"+<.4mm,.8mm>} **\crv{~*=<5pt>{.} "c"+(10,5) & "c"+(5,8)},
"e"*++!L{\scriptstyle X},
"xdown"+<0mm,-4mm> \ar _{\pi} @{>} {"x"+<0mm,10mm>}
\end{xy}
\end{center}
\end{wrapfigure}
for every stratum $F_R \subset \Xstar$.

Recall that for every $P\in\Pl$ with $P^\dag=R$ there is a factorization
$X_P=X_{P,\l}\times F_R$ and that $\pi|_{X_P}$ is simply projection onto
the second factor.  Thus $\pi^{-1}(x) = \coprod_{P^\dag=R} X_{P,\l}\times
\{x\} = \Xhat_{R,\l}\times \{x\}$ and we let $\ihat_{R,\l}: \Xhat_{R,\l}
\cong \pi^{-1}(x) \hookrightarrow \Xhat$ be the inclusion.  Since
$H^i(i_x^*\pi_*\M) = H^i(\Xhat_{R,\l};\ihat_{R,\l}^*\M)$ and
$H^i(i_x^!\pi_*\M) = H^i( \Xhat_{R,\l};\ihat_{R,\l}^!\M)$ we can use
Theorem~\ref{thmVanishing} to see these vanish for $i >
d(\ihat_{R,\l}^*\M)$ and $i< c(\ihat_{R,\l}^!\M)$ respectively.  Thus the
following theorem implies that \eqref{VanishingCovanishingCondition} holds
(and hence completes the proof of the Main Theorem):

\begin{thm}
\label{thmFunctoriality}
Let $\M$ be an $E$-micro-pure $\L$-module and let
$F_R$ be a stratum of the Baily-Borel-Satake compactification $\Xstar$.  Then
\begin{equation*}
d(\ihat_{R,\l}^*\M) \le \tfrac12\codim F_R -1 \qquad\text{and}\qquad
c(\ihat_{R,\l}^!\M) \ge \tfrac12\codim F_R +1 \ .
\end{equation*}
\end{thm}

The same result holds if $D$ is an equal-rank symmetric space and $\Xstar$
is a real equal-rank Satake compactification as in
\S\ref{sectGeneralization}.  This theorem is actually a special case of a
more general result on the functoriality of micro-support: for $\M$ an
arbitrary $\L$-module and $\Xstar$ a real equal-rank Satake
compactification as in \S\ref{sectGeneralization}, the theorem gives a
bound on $\mS(\ihat_{R,\l}^*\M)$ and $\mS(\ihat_{R,\l}^!\M)$ in terms of
$\mS(\M)$.

Since $\WC(\Xhat;E)$ is also $E$-micro-pure, the same argument yields a new
proof of the main result of \cite{refnGoreskyHarderMacPherson} (and in
fact a generalization to real equal-rank Satake compactifications).

\section{Example/application: ordinary cohomology}
\label{sectTilouine}
As another application of $\L$-modules we consider
the ordinary cohomology $H(X;\EE)$ or $H(\G;E)$ with coefficients in
$E\in\IrrRep(G)$.  This is the cohomology $H(\Xhat;\M)$ for the $\L$-module
$\M=i_{G*}E$ which has $E_G=E$ and $E_P=0$ for $P\neq G$ (see
\S\ref{sectExamples}(\ref{ExampleOne})).

We calculate the micro-support of $i_{G*}E$.  Since $i_Q^!i_{G*}E=E_Q$
we see that
\begin{equation*}
H(i_P^*\ihat_Q^!i_{G*}E) = \begin{cases} H(\n_P;E) & \text{for $Q=G$\ ,} \\
			      0		& \text{for $Q\neq G$\ .}
			\end{cases}
\end{equation*}
Thus for $V\in \IrrRep(L_P)$ to be in $\mS(i_{G*}E)$ it must be an
irreducible component of $H(\n_P;E)$ satisfying $(V|_{M_P})^* \cong
\overline { V|_{M_P}}$ and $(\xi_V+\r,\al)\le 0$ for all $\al\in \D_P$
(since $Q=G$ implies $Q_V'=G$).  The essential micro-support will consist
of such $V$ satisfying in addition the strict inequalities $(\xi_V+\r,\al)<
0$.

Let $\lambda$ be the highest weight of $E$.  As in
\S\ref{sectIHMicroPurity}, the irreducible components of $H(\n_P;E)$ are
the modules $V_{w(\lambda+\r)-\r}\in\IrrRep(L_P)$ with highest weight
$w(\lambda+\r)-\r$ for $w\in W_P$.  Let $\t_P:\h_\CC^*\to \h_\CC^*$
transform the highest weight of a representation of $L_P$ into the highest
weight of its complex conjugate contragredient; we assume that $\h =
\hb_{P} + \sa_P=\hb_{P,\k} + \hb_{P,\p} + \sa_P$ is a fundamental Cartan
subalgebra of $\levi_P$ equipped with a compatible order so that $\t_P$ is
simply the Cartan involution \cite{refnBorelCasselman}.  We can now
reexpress our calculation as
\begin{equation*}
\begin{split}
\emS(i_{G*}E)=\coprod_P \{\, V_{w(\lambda+\r)-\r} &\mid
\text{$w\in W_P$, $(w(\lambda+\r),\al)<0$ for all $\al\in\D_P$, and} \\
&\qquad \t_P(w(\lambda+\r)|_{\hb_P})=w(\lambda+\r)|_{\hb_P}\,\}\ .
\end{split}
\end{equation*}
(In the last equation we have used the fact that
$\t_P(\r|_{\hb_P})=\r|_{\hb_P}$.)  Furthermore since $V=
V_{w(\lambda+\r)-\r}$ occurs in $H(i_P^*\ihat_Q^!i_{G*}E)$ in degree $\l(w)$
we see that
\begin{equation}
\tilde c(V;i_{G*}E) = \tfrac12(\dim D_P -D_P(V))+ \l(w) \ .
\label{cVFormula}
\end{equation}

We use Lemma ~\ref{lemBasicLemma} to estimate $\l(w)$, however now we need
the term $\dim \n_P(V)$.  To define it, recall we have defined
$L_P(\u)\subseteq L_P$ in \S\ref{sectVanishingTheorem}\ to have roots
$\g\perp \u=w(\lambda+\r)-\r$.  Since $(w(\lambda+\r)-\r)|_{\hb_P}$ is
invariant under $\t_P$, the roots of $L_P(\u)$ are stable under $\t_P$.
Thus given an $L_P(\u)$-irreducible submodule of $\n_{P\CC}$, the transform
by $-\t_P$ of its weights are the weights of another $L_P(\u)$-irreducible
submodule of $\n_{P\CC}$.  Define $\n_P(\u)$ to be the sum of the
$L_P(\u)$-irreducible submodules of $\n_P$ whose weights are stable under
$-\t_P$.  Choose a compatible ordering for which $\dim \n_P(\u)$ is
maximized and let $\n_P(V)=\n_P(\u)$.  Note that $\n_P(V)$ contains the
root spaces of the positive $(-\t_P)$-invariant roots, that is, the real
roots.

We now make two assumptions: that $D$ is Hermitian symmetric, or more
generally equal-rank, and that $E$ has regular highest weight
$\lambda$.  By the first assumption the Lie algebra of $\lsp0
G(\RR)$ also possesses a compact Cartan subalgebra and therefore by
the Kostant-Sugiura theory of conjugacy classes of Cartan subalgebras
\cite{refnKostantConjugacyCartan}, \cite{refnSugiura},
\cite{refnSugiuraCorrection} there must exist at least $\dim
\hb_{P,\p} + \dim \sa_P - \dim \sa_G$ orthogonal real roots.  Thus
\begin{equation}
\dim \n_P(V)\ge \dim \hb_{P,\p} + \dim \sa_P - \dim \sa_G\ .
\label{nPVEstimate}
\end{equation}
On the other hand, note that if $\g^\vee=2\g/(\g,\g)$ then $(\r,\g^\vee)=1$
if and only if $\g$ is simple.  Consequently for $\g$ a simple root of
$L_P$ in any compatible ordering we have
\begin{equation*}
\begin{split}
\text{$\g$ is a root of $L_P(\u)$} &\Longleftrightarrow
(w(\lambda+\r),\g^\vee)=(\r,\g^\vee) \Longleftrightarrow
(\lambda+\r,w^{-1}\g^\vee)=1 \\
& \Longleftrightarrow (\lambda,w^{-1}\g)=0
\text{ and $w^{-1}\g$ is simple.}
\end{split}
\end{equation*}
Thus the second assumption implies that $L_P(\u)=H$, the Cartan
subgroup, and hence
\begin{equation}
\dim D_P(V) = \dim \hb_{P,\p}\ .
\label{DPVFormula}
\end{equation}

Lemma ~\ref{lemBasicLemma}(\ref{BasicLemmaPartOne}) and equations
\eqref{cVFormula}--\eqref{DPVFormula} yield the estimate $\tilde
c(V;i_{G*}E) \ge \tfrac12(\dim D_P +\dim \sa_P + \dim \n_P -\dim \sa_G) =
\tfrac12 \dim X$.  Thus Theorem ~\ref{thmVanishing} implies
\begin{thm}
\label{thmTilouine}
If $X$ is an arithmetic quotient of a Hermitian or
equal-rank symmetric space and $E$ has regular highest weight then
$H^i(X;\EE)=0$ for $i<\tfrac12\dim X$.
\end{thm}

This resolves a question posed by Tilouine during the Automorphic Forms
Semester.  For the case $G=R_{k/\QQ}\GSp(4)$ where $k$ is a totally real
number field the theorem is proven in \cite{refnTilouineUrban} using
results of Franke.  For applications of the theorem see \cite{refnMauger},
\cite{refnMokraneTilouine}.  While this paper was being prepared we heard
that Li and Schwermer also had a proof of the theorem.%
\footnote{Added Oct.~2003: See \cite{refnLiSchwermer}.  The methods are
  completely different.  They show vanishing in the range $i<\tfrac12(\dim
  X - (\rank \lsp0G - \rank K))$ without assuming $D$ is equal-rank.  This
  stregthened theorem also follows from the methods of the present paper:
  if $D$ is not equal-rank, equation \eqref{nPVEstimate} remains true
  provided we subtract $\rank \lsp0G - \rank K$ from the right-hand side.
  }

A vanishing range for the case where $E$ does not have regular highest
weight may be obtained by replacing \eqref{nPVEstimate} and
\eqref{DPVFormula} by the more subtle estimate on $\dim \n_P(V)$ given in
\cite{refnSaperLModules}.


\bibliography{ihp}

\newcommand{\SortNoop}[1]{}
\providecommand{\bysame}{\leavevmode ---\ }
\providecommand{\og}{``}
\providecommand{\fg}{''}
\providecommand{\smfandname}{et}
\providecommand{\smfedsname}{\'eds.}
\providecommand{\smfedname}{\'ed.}
\providecommand{\smfmastersthesisname}{M\'emoire}
\providecommand{\smfphdthesisname}{Th\`ese}
\begin{thebibliography}{10}

\bibitem{refnArthur}
{\scshape J.~Arthur} -- {\og The {$L^2$}-{L}efschetz numbers of {H}ecke
  operators\fg}, \emph{Invent. Math.} \textbf{97} (1989), p.~257--290.

\bibitem{refnBailyBorel}
{\scshape W.~Baily {\normalfont \smfandname} A.~Borel} -- {\og Compactification
  of arithmetic quotients of bounded symmetric domains\fg}, \emph{Ann. of
  Math.} \textbf{84} (1966), p.~442--528.

\bibitem{refnBorelCasselman}
{\scshape A.~Borel {\normalfont \smfandname} W.~Casselman} -- {\og
  {$L^2$}-cohomology of locally symmetric manifolds of finite volume\fg},
  \emph{Duke Math. J.} \textbf{50} (1983), p.~625--647.

\bibitem{refnBorelSerre}
{\scshape A.~Borel {\normalfont \smfandname} J.-P. Serre} -- {\og Corners and
  arithmetic groups\fg}, \emph{Comment. Math. Helv.} \textbf{48} (1973),
  p.~436--491.

\bibitem{refnCasselman}
{\scshape W.~Casselman} -- {\og {$L^2$}-cohomology for groups of real rank
  one\fg}, Representation Theory of Reductive Groups (P.~Trombi, \smfedname),
  Birkh{\"a}user, Boston, 1983, p.~69--84.

\bibitem{refnCasselmanGeometricRationality}
\bysame , {\og Geometric rationality of {S}atake compactifications\fg},
  Algebraic groups and {L}ie groups, Austral. Math. Soc. Lect. Ser.,
  \textbf{9}, Cambridge Univ. Press, Cambridge, 1997, p.~81--103.

\bibitem{refnvanEst}
{\scshape W.~T. van Est} -- {\og A generalization of the {C}artan-{L}eray
  sequence, {I}, {II}\fg}, \emph{Indag. Math.} \textbf{20} (1958), p.~399--413.

\bibitem{refnGoreskyHarderMacPherson}
{\scshape M.~Goresky, G.~Harder {\normalfont \smfandname} R.~MacPherson} --
  {\og Weighted cohomology\fg}, \emph{Invent. Math.} \textbf{116} (1994),
  p.~139--213.

\bibitem{refnGoreskyKottwitzMacPherson}
{\scshape M.~Goresky, R.~Kottwitz {\normalfont \smfandname} R.~MacPherson} --
  {\og Discrete series characters and the {L}efschetz formula for {H}ecke
  operators\fg}, \emph{Duke Math. J.} \textbf{89} (1997), p.~477--554.

\bibitem{refnGoreskyMacPhersonIHII}
{\scshape M.~Goresky {\normalfont \smfandname} R.~MacPherson} -- {\og
  Intersection homology {II}\fg}, \emph{Invent. Math.} \textbf{72} (1983),
  p.~77--129.

\bibitem{refnGoreskyMacPhersonWeighted}
\bysame , {\og Weighted cohomology of {S}atake compactifications\fg}, Centre de
  recherches math{\'e}matiques, preprint \#1593, 1988.

\bibitem{refnGoreskyMacPhersonTopologicalTraceFormula}
\bysame , {\og The topological trace formula\fg}, \emph{J. Reine Angew. Math.}
  \textbf{560} (2003), p.~77--150.

\bibitem{refnKashiwaraSchapira}
{\scshape M.~Kashiwara {\normalfont \smfandname} P.~Schapira} -- \emph{Sheaves
  on manifolds}, Springer-Verlag, Berlin, 1990.

\bibitem{refnKostantConjugacyCartan}
{\scshape B.~Kostant} -- {\og On the conjugacy of real {C}artan subalgebras\@.
  {I}.\fg}, \emph{Proc. Nat. Acad. Sci. U. S. A.} \textbf{41} (1955),
  p.~967--970.

\bibitem{refnKostant}
\bysame , {\og {L}ie algebra cohomology and the generalized {B}orel-{W}eil
  theorem\fg}, \emph{Ann. of Math.} \textbf{74} (1961), p.~329--387.

\bibitem{refnLiSchwermer}
{\scshape J.-S. Li {\normalfont \smfandname} J.~Schwermer} -- {\og On the
  {E}isenstein cohomology of arithmetic groups\fg}, \emph{Duke Math. J.}
  \textbf{123} (2004), no.~1, p.~141--169.

\bibitem{refnLooijenga}
{\scshape E.~Looijenga} -- {\og {$L^2$}-cohomology of locally symmetric
  varieties\fg}, \emph{Compositio Math.} \textbf{67} (1988), p.~3--20.

\bibitem{refnMauger}
{\scshape D.~Mauger} -- {\og Alg\`ebres de {H}ecke quasi-ordinaires
  universelles\fg}, \emph{Ann. Sci. \'Ecole Norm. Sup. (4)} \textbf{37} (2004),
  no.~2, p.~171--222.

\bibitem{refnMokraneTilouine}
{\scshape A.~Mokrane {\normalfont \smfandname} J.~Tilouine} -- {\og Cohomology
  of {S}iegel varieties with {$p$}-adic integral coefficients and
  applications\fg}, \emph{Ast\'erisque} (2002), no.~280, p.~1--95, Cohomology
  of Siegel varieties.

\bibitem{refnRaghunathan}
{\scshape M.~S. Raghunathan} -- {\og Vanishing theorems for cohomology groups
  associated to discrete subgroups of semisimple {L}ie groups\fg}, \emph{Osaka
  J. Math.} \textbf{3} (1966), p.~243--256.

\bibitem{refnRaghunathanCorrection}
\bysame , {\og Corrections to ``{V}anishing theorems~\ldots''\fg}, \emph{Osaka
  J. Math.} \textbf{16} (1979), p.~295--299.

\bibitem{refnRapoportLetterBorel}
{\scshape M.~Rapoport} -- 1986, letter to A. Borel.

\bibitem{refnRapoportNote}
\bysame , {\og On the shape of the contribution of a fixed point on the
  boundary\@. the case of {$\QQ$}-rank 1\fg}, unpublished, 1989.

\bibitem{refnRapoport}
\bysame , {\og On the shape of the contribution of a fixed point on the
  boundary: The case of {$\QQ$}-rank 1 \textup(with an appendix by {L.} {S}aper
  and {M.} {S}tern\textup)\fg}, The Zeta functions of {P}icard modular surfaces
  (R.~P. Langlands {\normalfont \smfandname} D.~Ramakrishnan, \smfedsname), Les
  Publications CRM, Montr{\'e}al, 1992, p.~479--488.

\bibitem{refnSaperTilings}
{\scshape L.~Saper} -- {\og Tilings and finite energy retractions of locally
  symmetric spaces\fg}, \emph{Comment. Math. Helv.} \textbf{72} (1997),
  p.~167--202.

\bibitem{refnSaperLModules}
\bysame , {\og {$\L$}-modules and micro-support\fg}, \emph{Ann. of Math.},
  accepted subject to revision, {\tt math.RT/0112251}, 2001.

\bibitem{refnSaperGeometricRationality}
\bysame , {\og Geometric rationality of equal-rank {S}atake
  compactifications\fg}, \emph{Math. Res. Lett.} \textbf{11} (2004),
  p.~653--671.

\bibitem{refnSaperSternTwo}
{\scshape L.~Saper {\normalfont \smfandname} M.~Stern} -- {\og
  {$L_2$}-cohomology of arithmetic varieties\fg}, \emph{Ann. of Math.}
  \textbf{132} (1990), p.~1--69.

\bibitem{refnSatakeCompact}
{\scshape I.~Satake} -- {\og {\SortNoop{1}}{O}n representations and
  compactifications of symmetric {R}iemannian spaces\fg}, \emph{Ann. of Math.}
  \textbf{71} (1960), p.~77--110.

\bibitem{refnSatakeQuotientCompact}
\bysame , {\og {\SortNoop{2}}{O}n compactifications of the quotient spaces for
  arithmetically defined discontinuous groups\fg}, \emph{Ann. of Math.}
  \textbf{72} (1960), p.~555--580.

\bibitem{refnSugiura}
{\scshape M.~Sugiura} -- {\og Conjugate classes of {C}artan subalgebras in real
  semisimple {L}ie algebras\fg}, \emph{J. Math. Soc. Japan} \textbf{11} (1959),
  p.~374--434.

\bibitem{refnSugiuraCorrection}
\bysame , {\og Correction to my paper: ``{C}onjugate classes of {C}artan
  subalgebras in real semisimple {L}ie algebras''\fg}, \emph{J. Math. Soc.
  Japan} \textbf{23} (1971), p.~379--383.

\bibitem{refnTilouineUrban}
{\scshape J.~Tilouine {\normalfont \smfandname} E.~Urban} -- {\og Several
  variable {$p$}-adic families of {S}iegel-{H}ilbert cusp eigensystems and
  their {G}alois representations\fg}, \emph{Ann. Sci. {\'E}cole Norm. Sup. (4)}
  \textbf{32} (1999), p.~499--574.

\bibitem{refnZuckerWarped}
{\scshape S.~Zucker} -- {\og {$L_2$} cohomology of warped products and
  arithmetic groups\fg}, \emph{Inv. Math.} \textbf{70} (1982), p.~169--218.

\bibitem{refnZuckerSatakeCompactifications}
\bysame , {\og {S}atake compactifications\fg}, \emph{Comment. Math. Helv.}
  \textbf{58} (1983), p.~312--343.

\bibitem{refnZuckerLtwoIHTwo}
\bysame , {\og {$L_2$}-cohomology and intersection homology of locally
  symmetric varieties, {II}\fg}, \emph{Compositio Math.} \textbf{59} (1986),
  p.~339--398.

\end{thebibliography}
\bibliographystyle{smfplain}
\end{document}